\documentclass[11pt]{article}
\usepackage{amsmath,amsfonts,amssymb}
\usepackage{float}
\usepackage{url}
\usepackage[english]{babel}
\usepackage{fullpage}

\usepackage{graphicx}
\newcommand{\setof}[1]{\left\{ {#1} \right\}}
\newcommand{\derp}[2]{\frac{\partial{#1}}{\partial{#2}}}
\newcommand{\derd}[2]{\frac{{\partial}^2{#1}}{\partial{#2}^2}}
\newcommand{\prodscal}[2]{\left( {#1} \mid {#2} \right)}

\def\diver{\mathop{\rm div}\nolimits}
\def\IR{\mathbb{R}}

\begin{document}
\title{Wind driven 3D Navier-Stokes circulation in the Atlantic}
\author{Olivier Besson, \\
Institut de math\'ematiques. \medskip\\
Julien Straubhaar, \\
Institut de math\'ematique and Centre d'hydrog\'eologie. \medskip \\
Universit\'e de Neuch\^atel, \\
11, rue Emile Argand, \\
CH-2000 Neuch\^atel. \medskip \\
e-mail: olivier.besson@unine.ch
}
\date{July 2012}
\maketitle
\begin{abstract}
A finite element method for the numerical solution of the anisotropic 
Navier-Stokes equations in shallow domain is presented. This method 
take into account aspect ratio in the hydrostatic approximation of the Navier-Stokes 
equations \cite{beslay,azthese,azgui}. A projection method 
\cite{guermond,shen} is used for the time discretization. The linear 
systems are solved via a some preconditioned conjugate algorithm, 
well adapted to massively parallel computers 
\cite{julien_these,julien_precond,julien_precond_paral}. Some 
results are presented for the wind driven water circulation in the
North Atlantic.
\end{abstract}
\section{Introduction}
\label{sec:intro}

The numerical simulation of wind-driven currents in the ocean is a common and current subject. Most
of the used models are based on the hydrostatic approximation and anisotropic viscosity. 
A Galerkin method for the primitive equations in shallow domains is presented in \cite{lw1,lw2}.
Prismatic finite elements are used, and numerical example are presented for the English Channel. 
In \cite{gwl}, a baroclinic circulation model for the North Atlantic is presented.
Spherical-polar coordinated are used . 
A study of the interaction between a baroclinic and a barotropic model is given in \cite{oe}. It is
based on the planetary geostrophic equations. 
Finally in \cite{win} a barotropic model is used for the study of the wind-driven circulation in an
elongated basin.
Most of the anisotropic models used in oceanography can be mathematically justified; see e.g.
\cite{beslay,azthese,tezi,azgui}.

Since the publication of the M.S. Lozier paper: Deconstructing the Conveyor
Belt \cite{lozier}, it is of importance to  have a good description of the
wind-driven circulation in oceans. This paper is devoted to the numerical
simulation of such a circulation in the North Atlantic, with a particular
attention to the 3D aspects; including down and up-welling.

In this paper, the full system of Navier-Stokes equations is considered; it is organized as
follows. 
In section \ref{sec:ns_ani}, we recall the setting of the anisotropic Navier-Stokes equations in the
context of thin domains. 
Then a description of a projection method is presented in section \ref{sec:proj}. 
A weak formulation and a finite element discretization are given in sections \ref {sec:weak} and
\ref{sec:fem}. 
We end this paper with section \ref{sec:atl} where some numerical simulations obtained with
our method are presented. Some wind driven 3D Navier-Stokes circulation in the North
Atlantic is reported.
\section{Anisotropic Navier-Stokes equations}
\label{sec:ns_ani}
Water flows in oceanography and
limnology are governed by the Navier-Stokes equations. For some numerical
simulations, asymptotic models are in current use (see \cite{pedlo,zeytou}).
These are all based on the fact that
the horizontal dimensions of the considered domain are much 
larger than the vertical one. 
Let $d$ be the horizontal width, and $h$ be the depth.
The simplest model using the fact that 
\[\epsilon = \frac{h}{d}\]
is very small is the hydrostatic model.
In this model, we take care of turbulence effects by setting an 
anisotropic viscosity, much smaller in the vertical direction than in 
the horizontal one (see \cite{pedlo}).

\begin{figure}[H]
\begin{center}
\includegraphics[width=12cm]{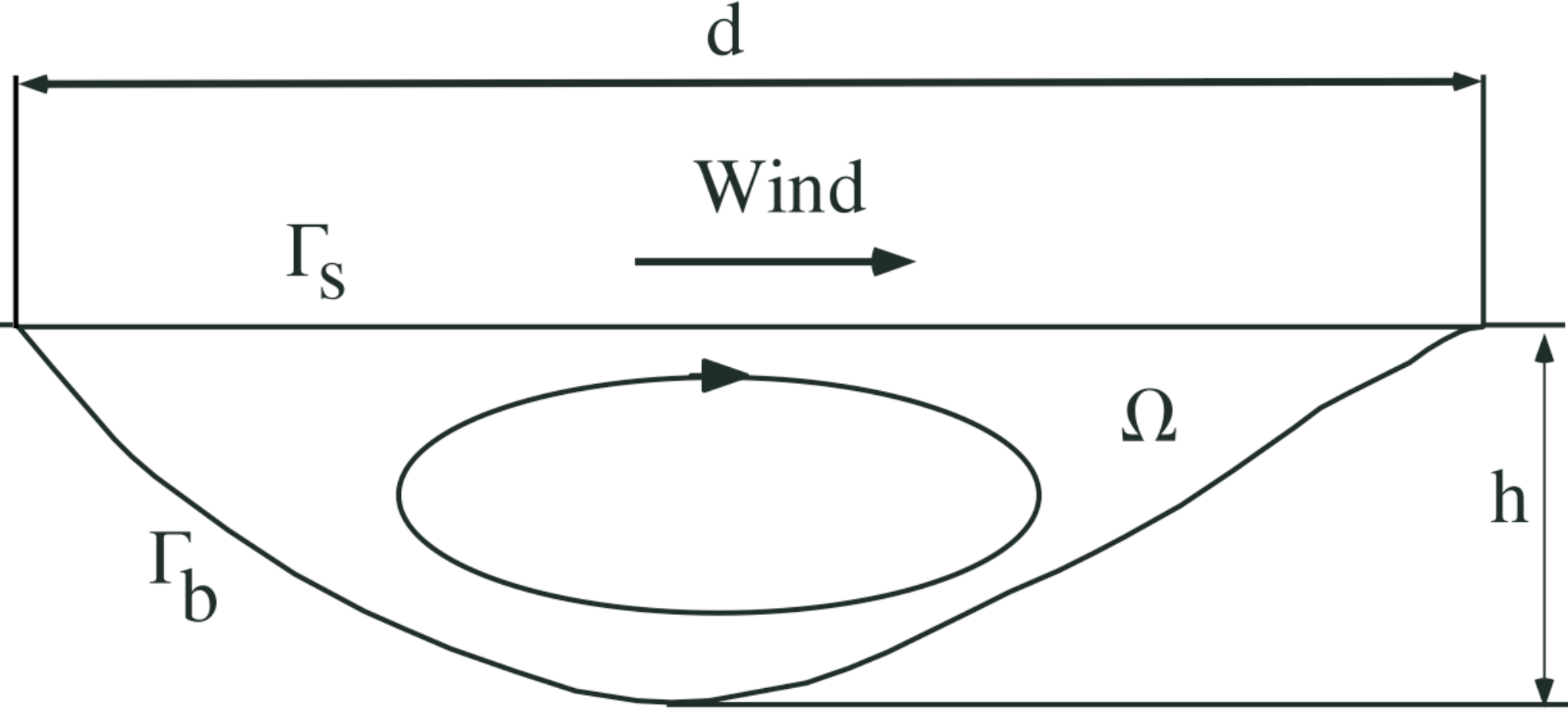}
\end{center}
\caption{Illustration of the shallow domain}
\label{schema}
\end{figure}

Let $ \Omega = \Omega_\epsilon \subset \IR^{3}$ be the
domain defined by
\[\Omega = \setof{x = (x_1,x_2,x_3) \in \IR^{3},
 \; (x_1, x_2) \in \Gamma_s, 
\; - h(x_1,x_2) < x_3 < 0}\]
where $\Gamma_s$ is the surface of the domain 
and $h: \Gamma_s \rightarrow \IR$ is its depth. 
The bottom of the domain is defined by  $\Gamma_b = \partial \Omega \setminus \Gamma_s$ 
(Fig. \ref{schema}). 
It is assumed that the water motion is generated by horizontal tension,  
induced by some wind on the surface $\Gamma_s$. This motion is driven by 
the following anisotropic Navier-Stokes equations and is
influenced by the Coriolis force:
\begin{eqnarray}
u_t + \prodscal{u}{\nabla}u &=&
 \Delta_{\nu}u - 2\omega \wedge u - \nabla p  
 \quad {\rm in \; \Omega}, \; t>0,
 \label{prim_velo}\\
\diver\, u &=& 0 \quad {\rm in} \; \Omega, \; t>0,
\label{prim_div}
\end{eqnarray}
with the following boundary and initial conditions.
The bottom $\Gamma_b$ is subdivided into two parts $\Gamma_0$ with strictly
positive measure, and $\Gamma_1 = \Gamma_b \setminus \Gamma_0$, and
\begin{equation}
u = 0 \quad {\rm on} \;  \Gamma_0, \; t>0,
\label{prim_bot_0}
\end{equation}
\begin{equation}
\prodscal{u}{n} = 0, \; \prodscal{\sigma \cdot n}{\tau} = 0 \quad {\rm on} \; 
\Gamma_1, \; t>0,
\label{prim_bot_1}
\end{equation}
where $\sigma = \sigma_{ij}$ is the stress tensor, with $\sigma = -pI + \nabla
u \cdot \nu + \nu \cdot \nabla u^t$, $n$ is the unit outward normal vector to
the boundary, and $\tau$ is any tangential vector. Note that the part $\Gamma_1$ 
of the bottom corresponds to artificial limits of the domain.

On the surface $\Gamma_s$, tension conditions are considered
\begin{equation}
\nu_3 \derp{u_1}{x_3} = \theta_1, \quad  \nu_3 \derp{u_2}{x_3} = 
\theta_2, \quad  u_3 = 0 \quad  {\rm on} \; \Gamma_s, \; t>0.
\label{prim_trac}
\end{equation}
Finally the initial condition is
\begin{equation}
u(\cdot ,\,t=0)   =  0 \quad  {\rm in} \; \Omega.
\label{prim_ci}
\end{equation}
The following notations have been used:
$u = (u_1, u_2, u_3)$ is the fluid velocity,
$\omega = (0, 0, \omega_{3})$ is the angular velocity of the Earth 
(projected onto the vertical in local coordinates),
$p$ is the pressure,
$\nu = \mathrm{diag}(\nu_1, \nu_2, \nu_3)$ is the turbulent viscosity diagonal
tensor,
$\theta_1$ and $\theta_2$ are the tensions induced by the wind
and 
\[ \Delta_{\nu} \varphi = \sum_{j=1}^{3} \; \nu_j \; 
\derd{\varphi}{x_j}.\] 
\medskip
Let us do the following change of variables and functions
\[ x_1 := x_1, \quad   x_2 := x_2, \quad  x_3 := x_3/\epsilon \]
\[ v_1 := u_1, \quad  v_2 := u_2, \quad  v_3 := u_3/\epsilon \]
%
With this scale change we get
\begin{eqnarray*}
\derp{v_1}{t} + \prodscal{v}{\nabla v_1} - \Delta_{{\nu}^{\epsilon}}v_1 - f v_2
+ \derp{p}{x_1} &=& 0 \quad {\rm in} \; \Omega, \; t>0, \\
\derp{v_2}{t} + \prodscal{v}{\nabla v_2} - \Delta_{{\nu}^{\epsilon}}v_2 + f v_1
+ \derp{p}{x_2} &=&   0 \quad {\rm in} \; \Omega, \; t>0, \\
{\epsilon}^2 \;\{\derp{v_3}{t} + \prodscal{v}{\nabla v_3} - 
\Delta_{{\nu}^{\epsilon}}v_3\} + \derp{p}{x_3}  &=&  0 \quad 
{\rm in} \; \Omega, \; t>0, \\
\diver \,v &=&   0 \quad {\rm in} \; \Omega, \; t>0,
\end{eqnarray*}
\begin{equation*}
{\nu}_3 \derp{v_1}{x_3}  =  \epsilon \theta_1, \quad \nu_3 
\derp{v_2}{x_3} = \epsilon \theta_2, \quad  v_3  =  0     \quad  
{\rm on} \; {\Gamma}_s,  \; t>0,
\end{equation*}
\begin{equation*}
v = 0 \quad {\rm on} \;  \Gamma_0, \; t>0,
\end{equation*}
\begin{equation*}
\prodscal{v}{n} = 0, \; \prodscal{\sigma_{\epsilon} \cdot n}{\tau} = 0 
\quad {\rm on} \; \Gamma_1, \; t>0,
\end{equation*}
with ${\nu}^{\epsilon} = \mathrm{diag}({\nu}_1, {\nu}_2, {\nu_3}/{\epsilon^2})
$, 
and $f = 2{\omega}_3$. \\
Set 
\[
\nu_1 = \lambda_1, \; \nu_2 = \lambda_2, \; \nu_3 = \epsilon^{2}\,\lambda_3,
\]
\[\Theta_i = \epsilon^{-1}\,\theta_i \;\;\; i = 1, 2,\] 
\[ \lambda = (\lambda_1, \lambda_2, \lambda_3 ). \]
When $\epsilon \rightarrow 0$,  this problem becomes 
\emph{the hydrostatic approximation of Navier-Stokes equations.}
\begin{eqnarray}
\derp{v_1}{t} + \prodscal{v}{\nabla v_1} - \Delta_{{\lambda}}v_1 - f v_2 +
\derp{p}{x_1} &=&  0 \quad {\rm in} \; \Omega,  \; t>0,
\label{hydr_1}\\
\derp{v_2}{t} + \prodscal{v}{\nabla v_2} - \Delta_{\lambda}v_2 + f v_1 +
\derp{p}{x_2}  &=&   0 \quad {\rm in} \; \Omega,  \; t>0,
\label{hydr_2}\\
\derp{p}{x_3}   &=&  0 \quad {\rm in} \; \Omega,  \; t>0,
\label{hyrd_3}\\
\diver v  &=&   0 \quad {\rm in} \; \Omega,  \; t>0.
\label{hydr_div}
\end{eqnarray}
with
\begin{equation}
v_1 = v_2 = v_3 \cdot n_3  =  0 \quad {\rm on} \;  {\Gamma}_0,  \; t>0,
\label{hydr_bot0}
\end{equation}
\begin{equation}
\prodscal{v_H}{n_H} = 0, \; \prodscal{\sigma_H \cdot n_H}{\tau_H} = 0 
\quad {\rm on} \; \Gamma_1, \; t>0,
\label{hydr_bot1}
\end{equation}
where $w_H$ denotes the horizontal components of $w$,
\begin{equation}
{\lambda}_3 \derp{v_1}{x_3}  =   \Theta_1, \quad \lambda_3 
\derp{v_2}{x_3} =  \Theta_2, \quad  v_3  =  0     \quad  {\rm on} \; 
{\Gamma}_s ,  \; t>0,
\label{hydr_trac}
\end{equation}
\begin{equation}
v_{1}(\cdot ,\,t=0) = v_{2}(\cdot ,\,t=0)   =  0 \quad  {\rm in} \; \Omega.
\label{hydr_ci}
\end{equation}
This development shows that it is natural to use an anisotropic viscosity 
diagonal tensor with the third component of the order of $\epsilon^2$, compared 
with the others.
\section{Projection method}
\label{sec:proj}
In order to solve equations \eqref{prim_velo}-\eqref{prim_ci}, a 
velocity-correction projection method \cite{guermond,shen,guermond_shen} 
is used. Let us recall this method in our case. 
Set $u^0 = u(0), p^0 = p(0)$, and let $u^1$, and
$p^1$ be approximations of $u(\delta t)$ and $p(\delta t)$.
A BDF2 scheme is used for the time discretization.
 For $k \geq 1$, we look for $\tilde{u}^{k+1},
p^{k+1}$, and $u^{k+1}$ such that the prediction
$\tilde{u}^{k+1}$ of the velocity is a
solution of
\begin{eqnarray}
\frac{1}{\delta t}\left(\frac{3}{2} \tilde{u}^{k+1} - 2 u^k +
\frac{1}{2} u^{k-1}\right)-\Delta_\nu \tilde{u}^{k+1} & = &
-\prodscal{u^k}{\nabla}u^k-2\omega \wedge u^k
-\nabla p^{k} \; \text{ in }\Omega, \label{pred_velo}\\
\tilde{u}^{k+1} & = &  0 \; \text{ on }  \Gamma_0,\label{pred_cl_0} \\
\prodscal{\tilde{u}^{k+1}}{n} & = &  0 \; \text{ on }  \Gamma_1,\label{pred_cl_1_u} \\
\nu_3 \frac{\partial \tilde{u}_1^{k+1}}{\partial x_3} &=& \theta_1^{k+1}, \; \text{ on } \Gamma_s
\label{pred_tr_1} \\
\nu_3 \frac{\partial \tilde{u}_2^{k+1}}{\partial x_3} &=& \theta_2^{k+1}, \; \text{ on } \Gamma_s
\label{predc_tr_2} \\
\tilde{u}_3^{k+1} &=& 0 \; \text{ on } \Gamma_s, \label{pred_tr_3}
\end{eqnarray}
then the velocity $u^{k+1}$, and the pressure $p^{k+1}$ satisfy
\begin{eqnarray}
\frac{3}{2\delta t}\left(u^{k+1}-\tilde{u}^{k+1}\right) + \nabla \left(p^{k+1} -
p^{k}\right) &=& 0 \; \text{ in }\Omega, \label{correc_velo_u} \\
\diver u^{k+1} &=& 0 \; \text{ in }\Omega,\label{correc_inc_u} \\
\prodscal{u^{k+1}}{n} &=& 0 \; \text{ on } \partial \Omega. \label{correc_cl_u} 
\end{eqnarray}
If $q^{k+1} =p^{k+1}-p^{k}$, from equations \eqref{correc_velo_u} and \eqref{correc_inc_u} we have
\begin{eqnarray}
- \Delta q^{k+1} &=& -\frac{3}{2\delta t} \diver \tilde{u}^{k+1} \; \text{ in }\Omega, 
\label{press_cor} \\
\derp{q^{k+1}}{n} & = & 0  \; \text{ on } \partial \Omega. \label{press_cor_cl} 
\end{eqnarray}
If $q^{k+1}$ is known, then $p^{k+1} = q^{k+1}+p^{k}$, and 
$\displaystyle u^{k+1} = \tilde{u}^{k+1} + \delta t \frac{2}{3} \nabla q^{k+1}$.

Note that this method consists in computing a non divergence-free prediction $\tilde{u}^{k+1}$ of
the velocity, and then obtain $u^{k+1}$ as a divergence-free vector field via a pressure correction.
\section{Weak formulation}
\label{sec:weak}
A week formulation of problems \eqref{pred_velo}-\eqref{pred_tr_3} and
\eqref{press_cor}-\eqref{press_cor_cl} is the following.
Define
\begin{eqnarray*}
V &=& \setof{\varphi \in  H^1(\Omega); \; \varphi = 0 \; {\rm on } \; 
{\Gamma}_0}, \\
W &=& \setof{v \in V^3; \; 
\prodscal{v}{n} =0 \; {\rm on } \; \Gamma_s \cup \Gamma_1 }, \\
U &=& \setof{v \in H^1(\Omega)^3; \; 
\prodscal{v}{n} =0 \; {\rm on } \; \partial \Omega }.
\end{eqnarray*}
For $\theta_1, \, \theta_2 \in H^{-1/2}(\Gamma_s)$, if $k \geq 0$, assume that
$p^k \in H^1(\Omega)/\IR$, and $u^k, u^{k-1}$ in $U$ are given. 
We seek for $\tilde{u}^{k+1} \in W$ such that
\begin{multline}
\frac{3}{2\delta t}\int_{\Omega} \prodscal{\tilde{u}^{k+1}}{v}\, dx +
\int_{\Omega} \prodscal{\nu \cdot \nabla \tilde{u}^{k+1}}{\nabla v} \, dx = \\
\frac{1}{\delta t}\int_{\Omega} 
\left(2 \prodscal{u^k}{v} - \frac{1}{2}\prodscal{u^{k-1}}{v} \right) \, dx
- \int_{\Omega} \prodscal{\prodscal{u^k}{\nabla}u^k}{v} \, dx \\
- \int_{\Omega} 2\prodscal{\omega \wedge u^k}{v} \, dx
+ \int_{\Omega} p^k \, \diver v \, dx
+ \int_{\Gamma_s}\theta_1^{k+1} v_1\, ds + \int_{\Gamma_s}\theta_2^{k+1} v_2\, ds,
\label{fofa_velo_pred}
\end{multline}
for all $v \in W$. Then find $q^{k+1} \in H^1(\Omega)/\IR$ such that
\begin{equation}
\int_{\Omega} \prodscal{\nabla q^{k+1}}{\nabla \varphi} \, dx  = 
- \frac{3}{2\delta t}\int_{\Omega} \diver \tilde{u}^{k+1} \varphi \, dx,
\label{fofa_q_cor}
\end{equation}
for all $\varphi \in  H^1(\Omega)/\IR$. The pressure $p^{k+1} \in H^1(\Omega)/\IR$ is given by
\begin{equation}
\int_{\Omega} p^{k+1} \varphi \, dx = \int_{\Omega} (q^{k+1} + p^k) \varphi \, dx,
\label{fofa_pres_cor}
\end{equation}
for all $\varphi \in L^2(\Omega)$. Finally the velocity $u^{k+1}\in U$ is the solution of
\begin{equation}
\int_{\Omega} \prodscal{u^{k+1}}{v}\, dx = 
\int_{\Omega} \prodscal{\tilde{u}^{k+1}}{v} \, dx +
\delta t \frac{2}{3}\int_{\Omega} \prodscal{\nabla q^{k+1}}{v} \, dx,
\label{fofa_velo_cor}
\end{equation}
for all $v \in L^2(\Omega)^3$.
\section{Finite element discretization}
\label{sec:fem}
A finite element mesh $\tau_h$ of the domain $\Omega$ into 
brick elements is considered: $\overline{\Omega} = \cup_{K \in \tau_h} K$, where each geometric
elements $K$ is an hexahedron.
For $1 \leq j \leq 3$, set
\begin{eqnarray*}
V_{h} &=& \setof{\varphi_h \in  C(\overline{\Omega});
 \; \varphi_h|_K \in Q_2, \; \forall K \in \tau_h}, \\
W_h &=& \setof{v_h \in V_h^3; \; 
\prodscal{v_h}{n} =0 \; {\rm on } \, \Gamma_s \cup \Gamma_1; \;
v =0 \; {\rm on } \; \Gamma_0}, \\
U_h &=& \setof{v_h \in V_h^3; \; \prodscal{v_h}{n} =0 \; {\rm on } \, \partial \Omega}, \\
P_{h} &=& \setof{\varphi_h \in  C(\overline{\Omega});
\; \varphi_h|_K \in Q_1, \; \forall K \in \tau_h}.
\end{eqnarray*}
The following finite element discretization of the previous projection
method is used.

For $k \geq 0$, assume that
$p^k_h \in P_{h}$, and $u^k_h, u^{k-1}_h$ in $U_h$ are given. 
Then compute  $\tilde{u}^{k+1}_h \in W_h$ such that
\begin{multline}
\frac{3}{2\delta t}\int_{\Omega} \prodscal{\tilde{u}^{k+1}_h}{v_h}\, dx +
\int_{\Omega} \prodscal{\nu \cdot \nabla \tilde{u}^{k+1}_h}{\nabla v_h} \, dx = \\
\frac{1}{\delta t}\int_{\Omega} 
\left(2 \prodscal{u^k_h}{v_h} - \frac{1}{2}\prodscal{u^{k-1}_h}{v_h} \right) \, dx
- \int_{\Omega} \prodscal{\prodscal{u^k_h}{\nabla}u^k_h}{v_h} \, dx \\
- \int_{\Omega} 2\prodscal{\omega \wedge u^k_h}{v_h} \, dx
+ \int_{\Omega} p^k_h \, \diver v_h \, dx
+ \int_{\Gamma_s}\theta_1^{k+1} v_{h,1}\, ds + \int_{\Gamma_s}\theta_2^{k+1} v_{h,2}\, ds,
\label{fofa_velo_pred_h}
\end{multline}
for all $v_h \in W_h$. Then find $q_h^{k+1} \in P_h/\IR$ with
\begin{equation}
\int_{\Omega} \prodscal{\nabla q_h^{k+1}}{\nabla \varphi_h} \, dx  = 
- \frac{3}{2\delta t}\int_{\Omega} \diver \tilde{u}^{k+1}_h \varphi_h \, dx,
\label{fofa_q_cor_h}
\end{equation}
for all $\varphi_h \in  P_h$. The pressure $p_h^{k+1} \in P_h/\IR$ is given by
\begin{equation}
\int_{\Omega} p_k^{k+1} \varphi_h \, dx = \int_{\Omega} (q_h^{k+1} + p_h^k) \varphi_h \, dx,
\label{fofa_pres_cor_h}
\end{equation}
for all $\varphi_h \in P_h$. Finally the velocity $u_h^{k+1} \in U_h$ is the solution of
\begin{equation}
\int_{\Omega} \prodscal{u_h^{k+1}}{v_h}\, dx = 
\int_{\Omega} \prodscal{\tilde{u}^{k+1}_h}{v_h} \, dx +
\delta t \frac{2}{3}\int_{\Omega} \prodscal{\nabla q_h^{k+1}}{v_h} \, dx,
\label{fofa_velo_cor_h}
\end{equation}
for all $v \in V_h^3$.
\medskip \\
\textbf{Remark}
\begin{enumerate}
\item The solution of the linear systems \eqref{fofa_velo_pred_h}-\eqref{fofa_velo_cor_h} 
are performed via some preconditioned conjugate gradient methods. These methods are well adapted to
massively parallel computers.
\item For equation
\eqref{fofa_q_cor_h}, a coupled preconditionner: diagonal plus optimal conjugate
Grahm-Schmidt least squares preconditionner (DIAG + LS CGS OPT) is 
used, see \cite{julien_precond,julien_precond_paral,julien_these}. Moreover a Lagrange multiplier
is used in order to impose the zero mean of  $q_h^{k+1}$, see also \cite{bocleh}.
\item For the other
equations, an incomplete Cholesky IC0 preconditionner \cite{saad} is sufficient.
\end{enumerate}
\section{Application to the North Atlantic}
\label{sec:atl}

This section is devoted to the numerical simulation of the $3D$ water
circulation in the North Atlantic ocean, induced by some mean wind tensions. For
this, the methods presented in the previous sections are used. A parallel
software was developed to solve this kind of problems. 
The bathymetry, and the wind tensions where obtained from
the Mercator Ocean Project (http://www.mercator.fr). We are grateful to Jean-Marc Molines at
INP-Grenoble for his help.

The model into consideration does not take into account of the temperature, nor
the salinity. Our aim is, in a first step, to study the currents induced by the
wind tensions in the ocean. A particular attention is turned on the down, and upwelling near
the coast. Our model allows to obtain consistent results with reality. It demonstrate that the winds
are the main driving forces for the global dynamic in the oceans.

\begin{figure}[H]
\begin{center}
\includegraphics[width=12cm]{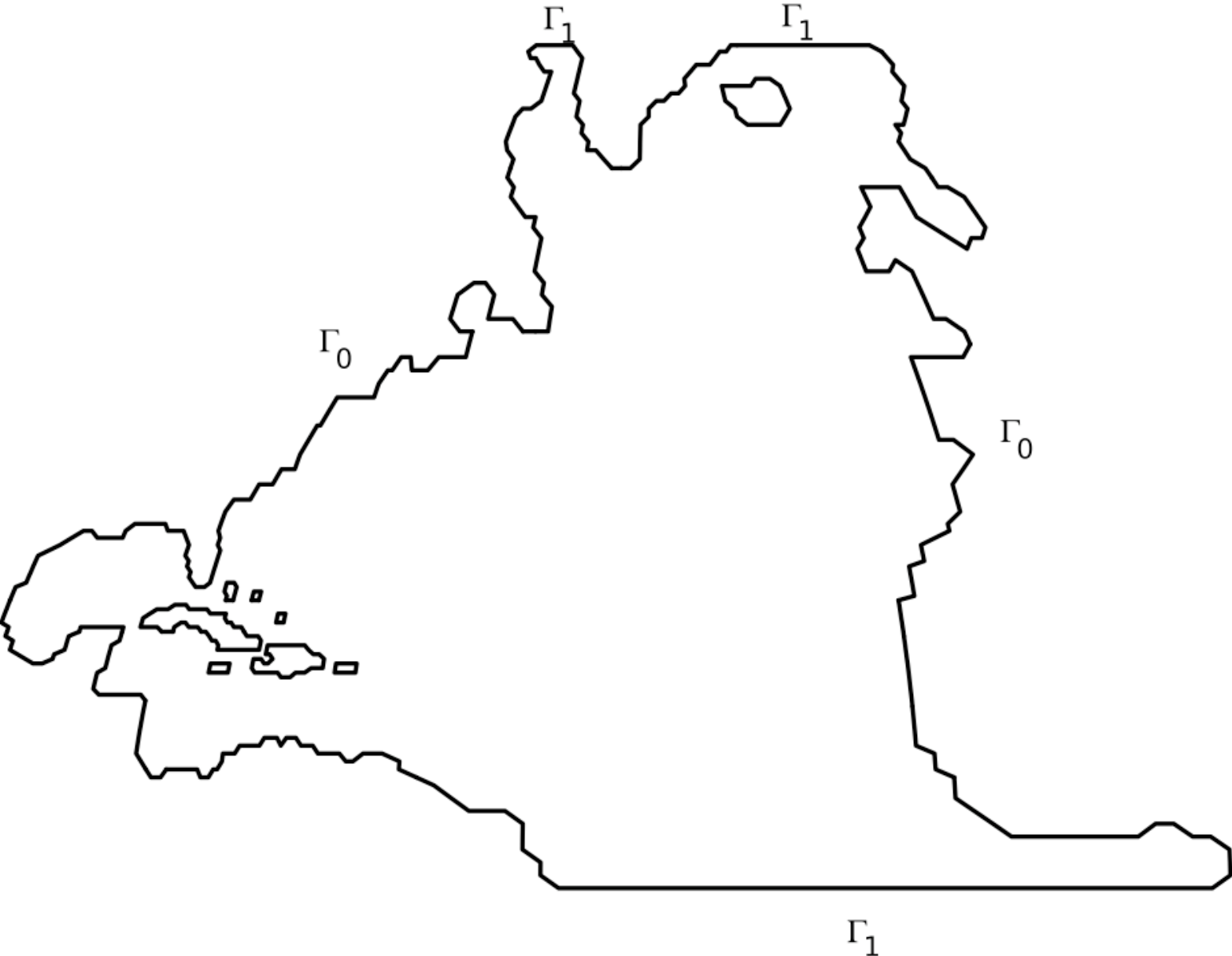}
\end{center}
\caption{Part of the Atlantic Ocean}
\label{atl}
\end{figure}

The part of the Atlantic ocean taken into account is delimited in the East by
the European and the African continents, in the West by the American continent,
in the South by the Equator, and in the North by the $70^{\circ}N$ parallel.
These two parts lead to some artificial boundary $\Gamma_1$ (see figure
\ref{atl}). On this part of the boundary, it is assumed that there is 
no tangential tensions, and that there is no water exchange with outside.

\begin{figure}[H]
\begin{center}
\includegraphics[width=12cm]{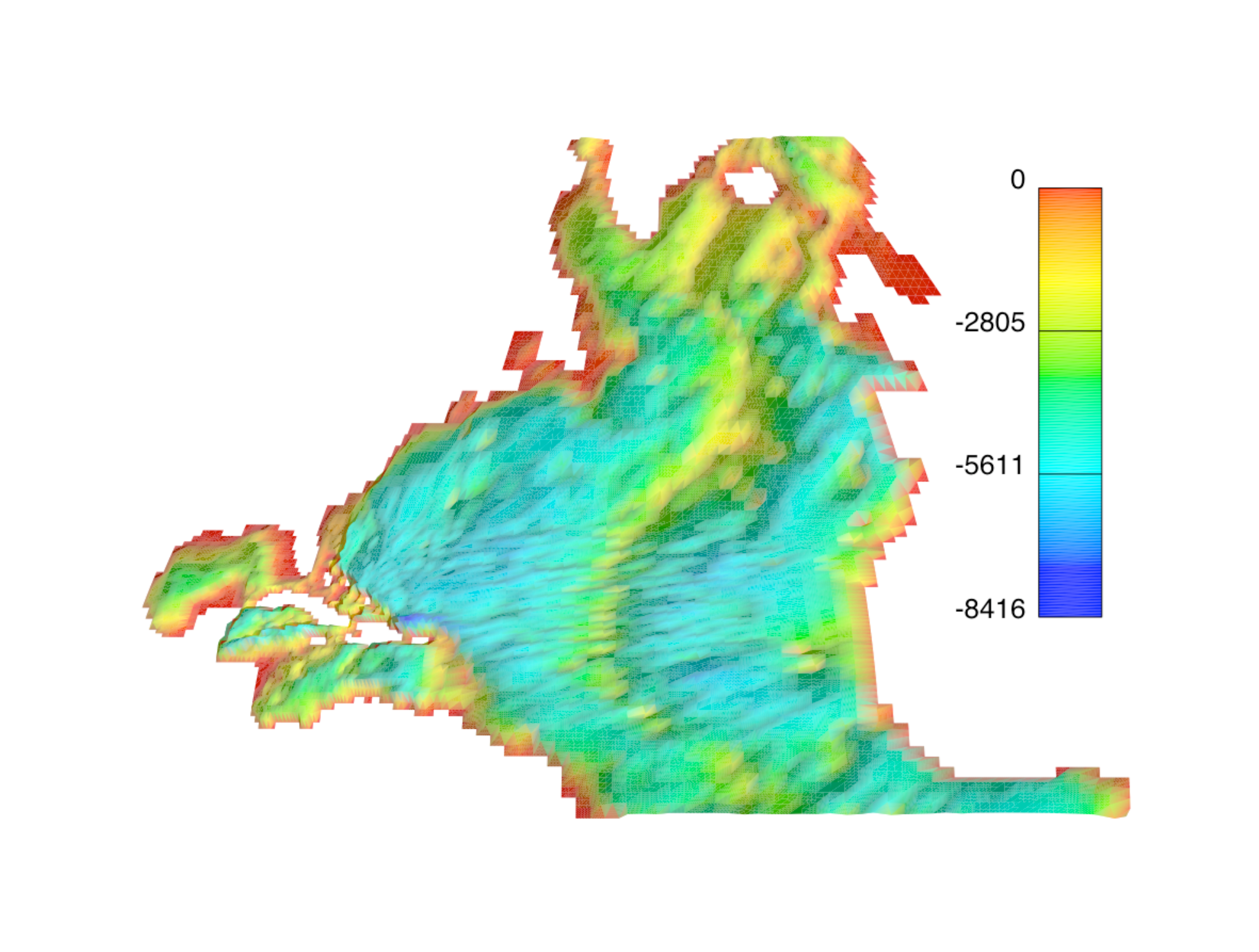}
\end{center}
\caption{Bathymetry of the Atlantic Ocean}
\label{bath_atl1}
\end{figure}

\begin{figure}[H]
\begin{center}
\includegraphics[width=12cm]{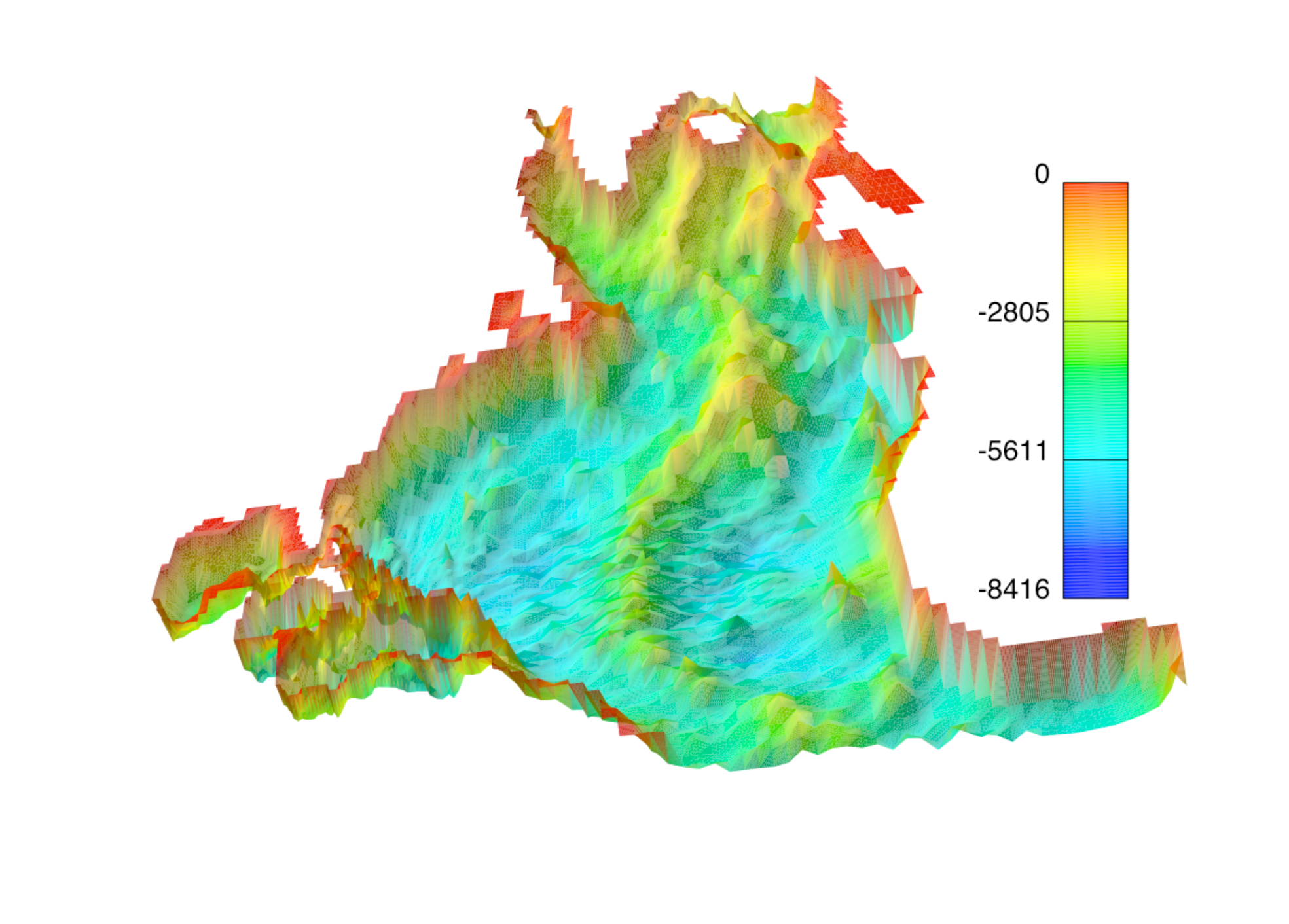}
\end{center}
\caption{Bathymetry of the Atlantic Ocean, with vertical scale}
\label{bath_atl2}
\end{figure}

\begin{figure}[H]
\begin{center}
\includegraphics[width=12cm]{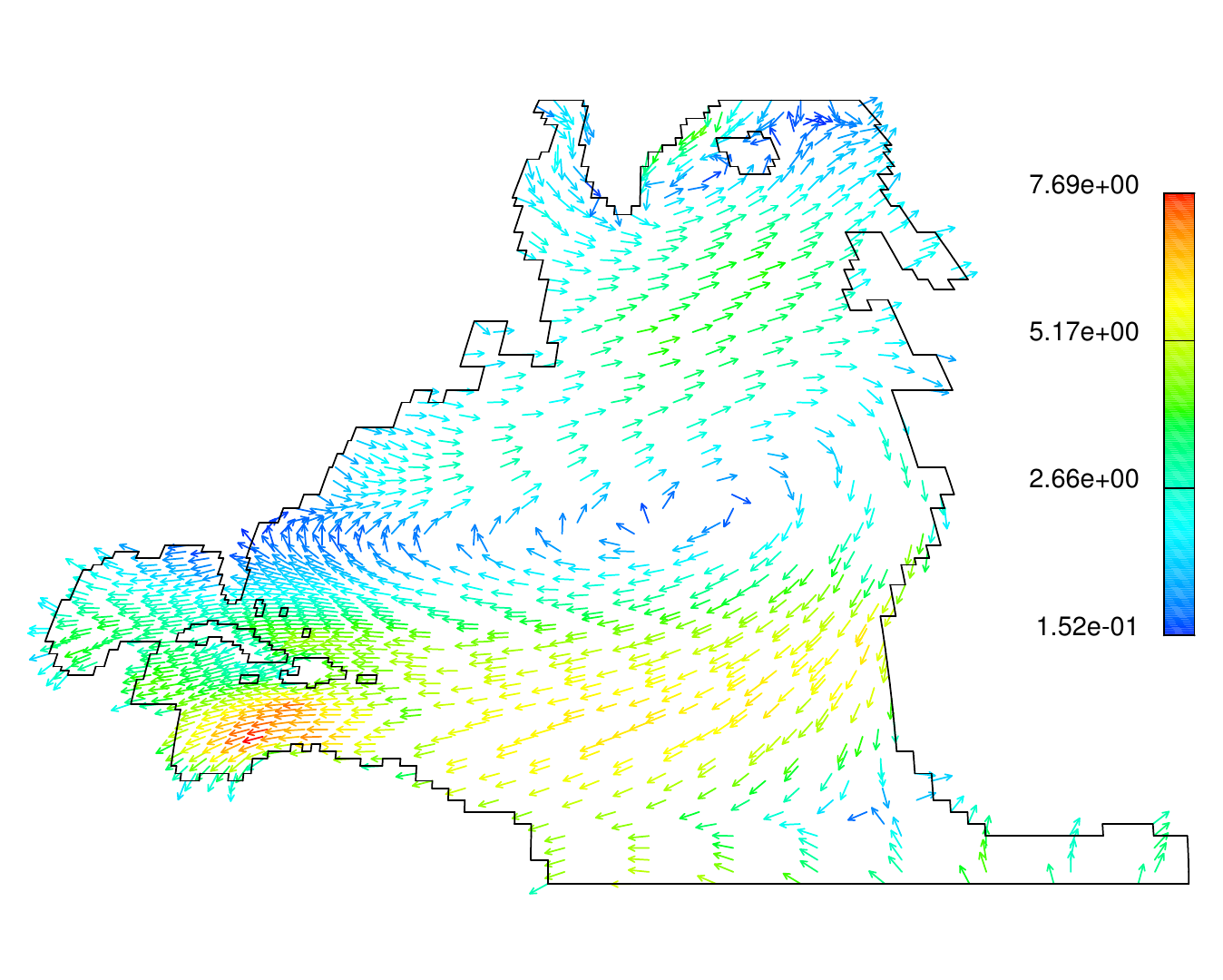}
\end{center}
\caption{Mean wind ERA15 (1979-1993)}
\label{wind}
\end{figure}

\begin{figure}[H]
\begin{center}
\includegraphics[width=12cm]{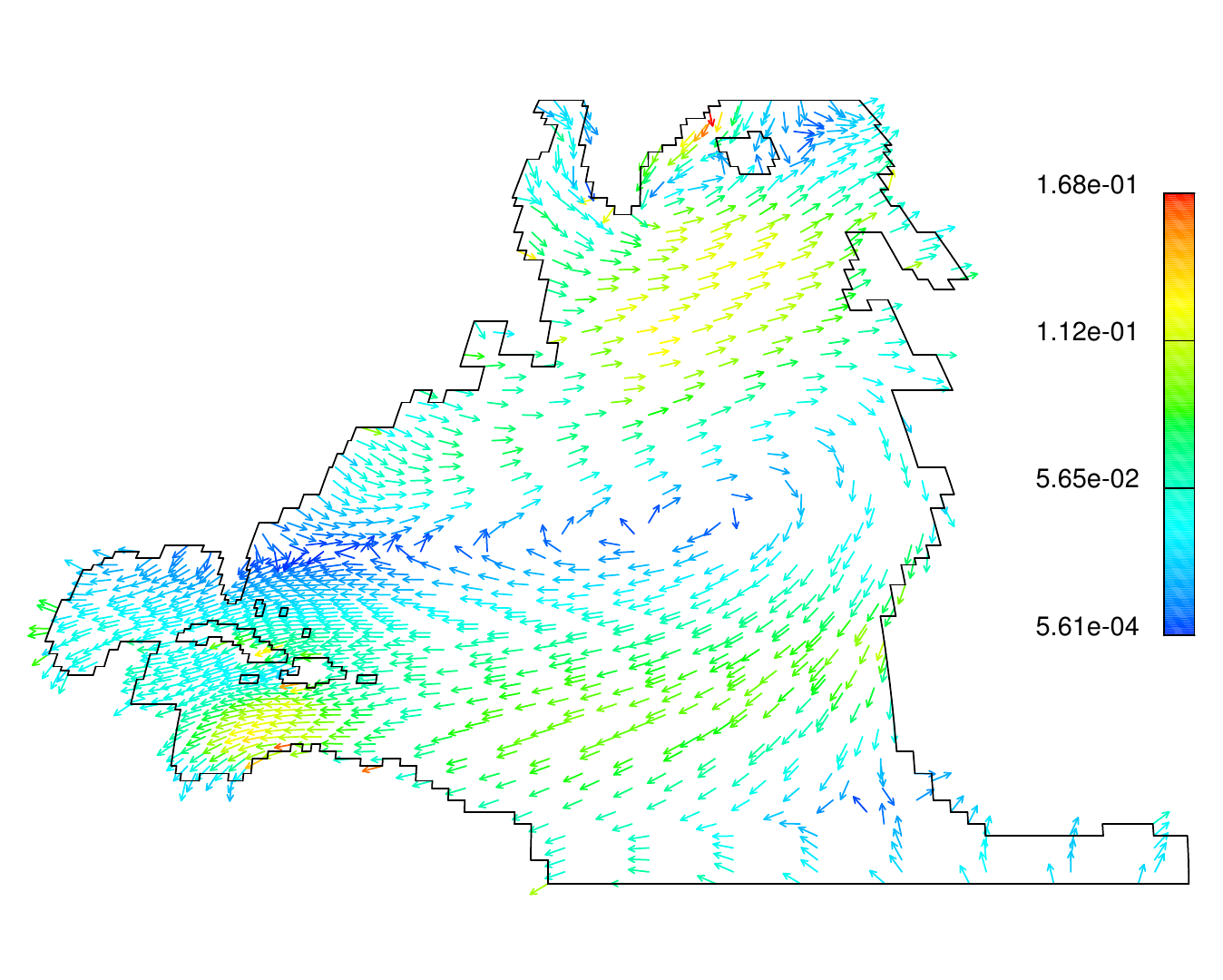}
\end{center}
\caption{Mean wind tensions ERA15 (1979-1993)}
\label{tension}
\end{figure}

The bathymetry is given in figure \ref{bath_atl1}, and \ref{bath_atl2}.
The mean wind during 15 years (1979-1993, ERA15) are represented in figure \ref{wind}
The mean wind tensions during this period are represented in figure \ref{tension}
The following physical constants are used in the equations:
\begin{itemize}
\item  The turbulent cinematic viscosity tensor is $\nu = (10^8, \, 10^8, \, 2.5 \cdot 10^2) \;
[m^2/s]$, it is based on the study of the Reynolds turbulent tensor.
\item  The time step is $\delta t = 2.592 \cdot 10^6 $ (= 1 month).
\item  The final time is set to 100 years
\end{itemize}

In the following figures, some numerical results for the 
currents at different depth are presented. 
\begin{figure}[H]
\begin{center}
\includegraphics[width=11cm]{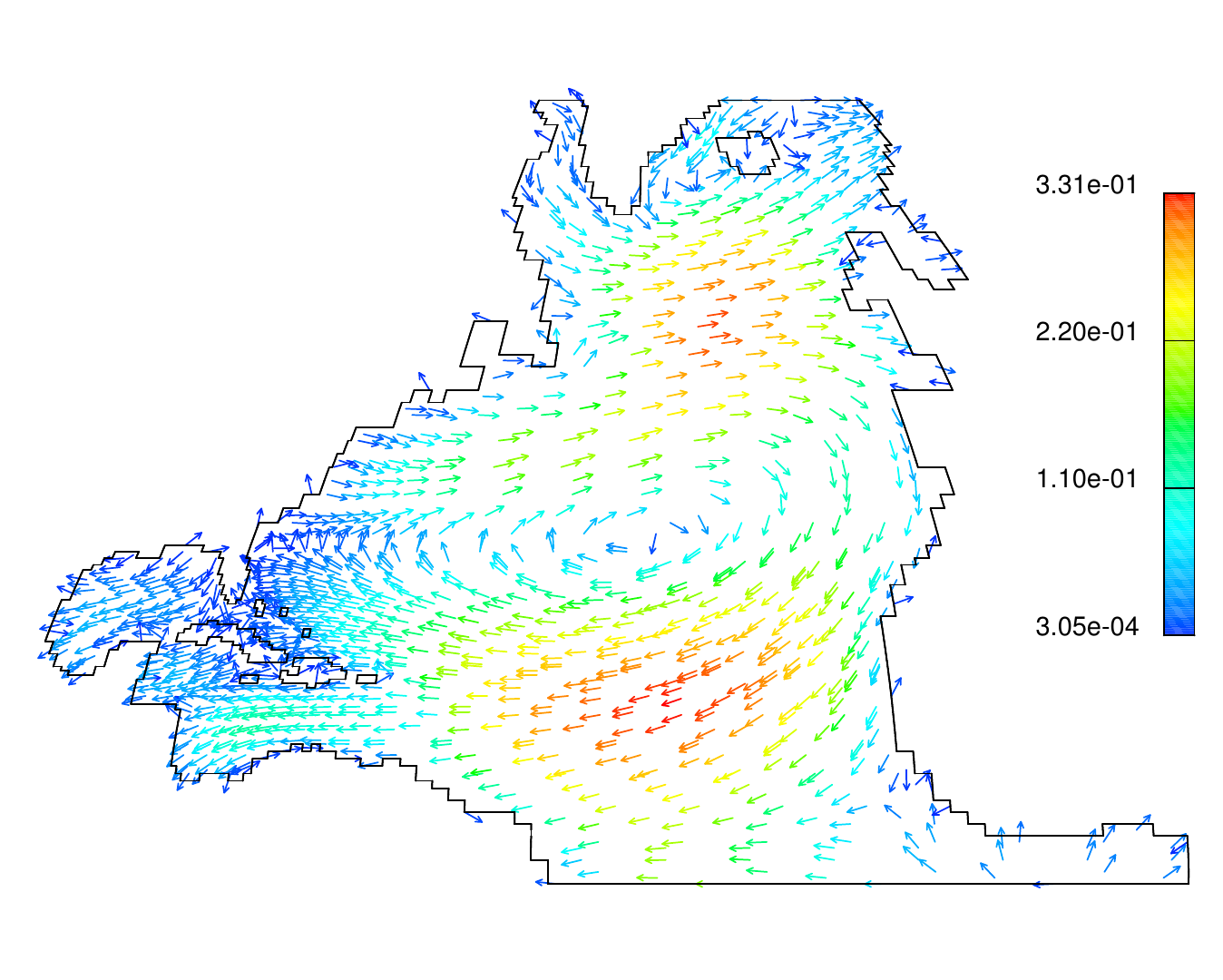}
\end{center}
\caption{Currents at the surface}
\label{velo0}
\end{figure}
\begin{figure}[H]
\begin{center}
\includegraphics[width=11cm]{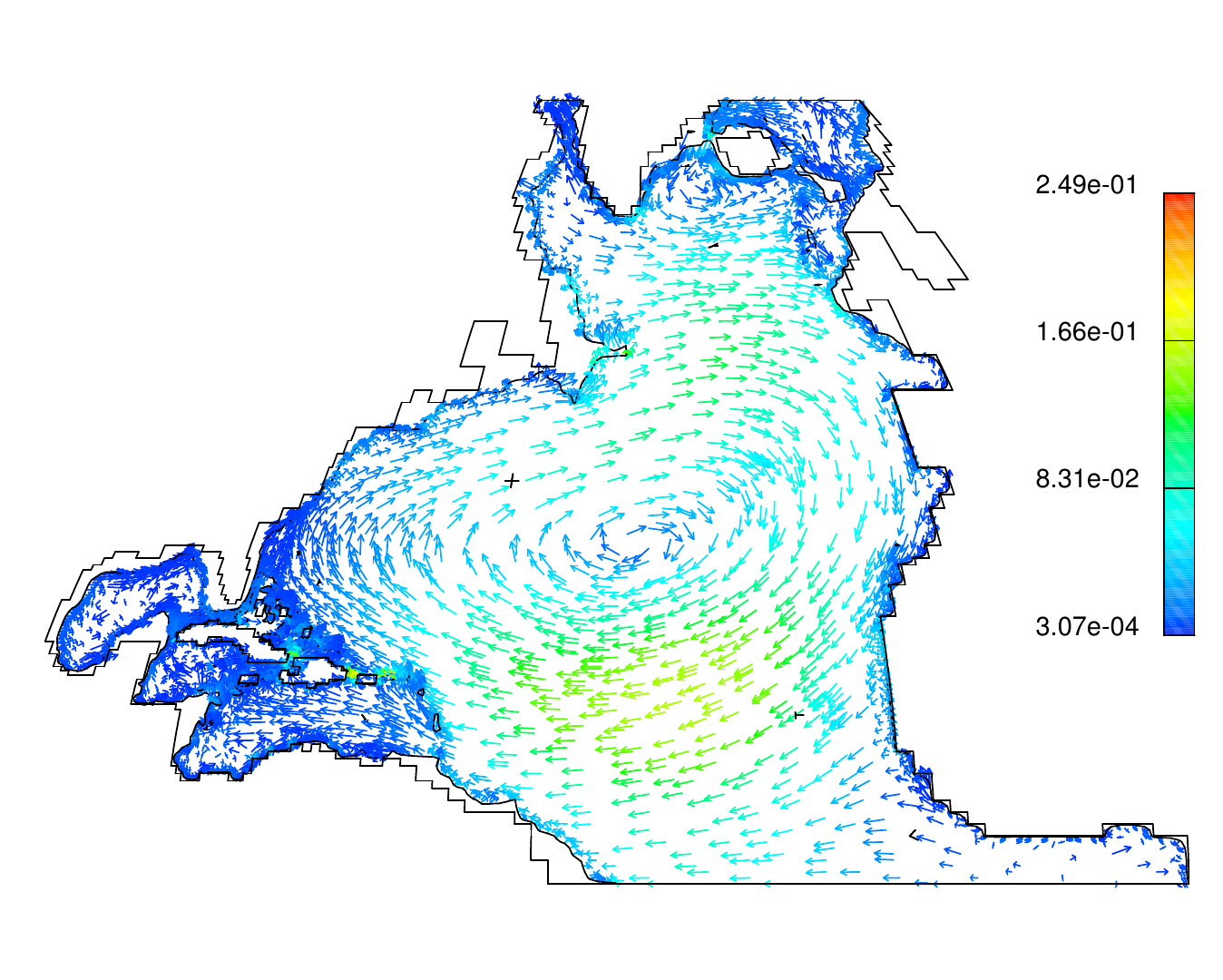}
\end{center}
\caption{Currents at 500 m depth}
\label{velo500}
\end{figure}
\begin{figure}[H]
\begin{center}
\includegraphics[width=11cm]{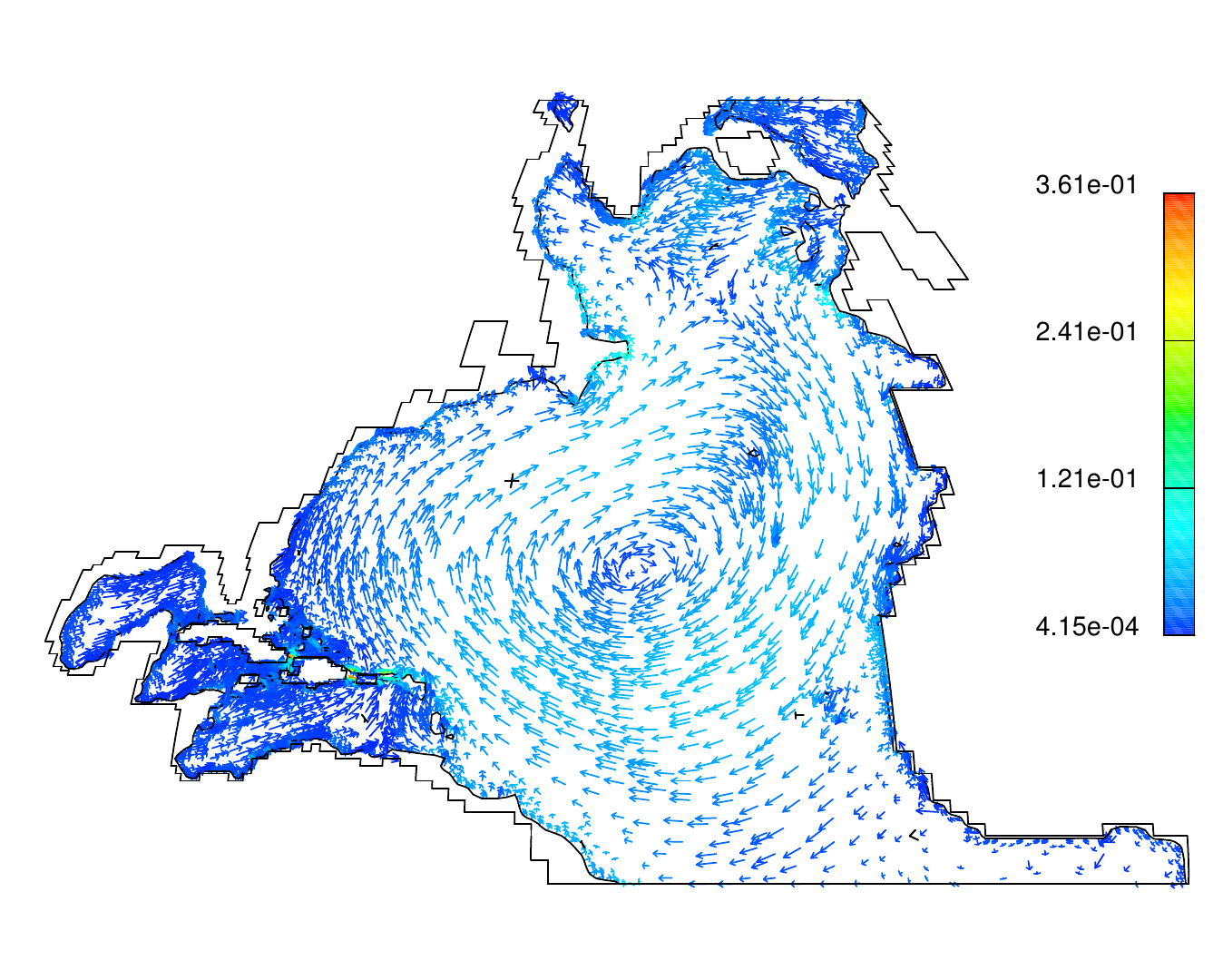}
\end{center}
\caption{Currents at 1000 m depth}
\label{velo1000}
\end{figure}
\begin{figure}[H]
\begin{center}
\includegraphics[width=11cm]{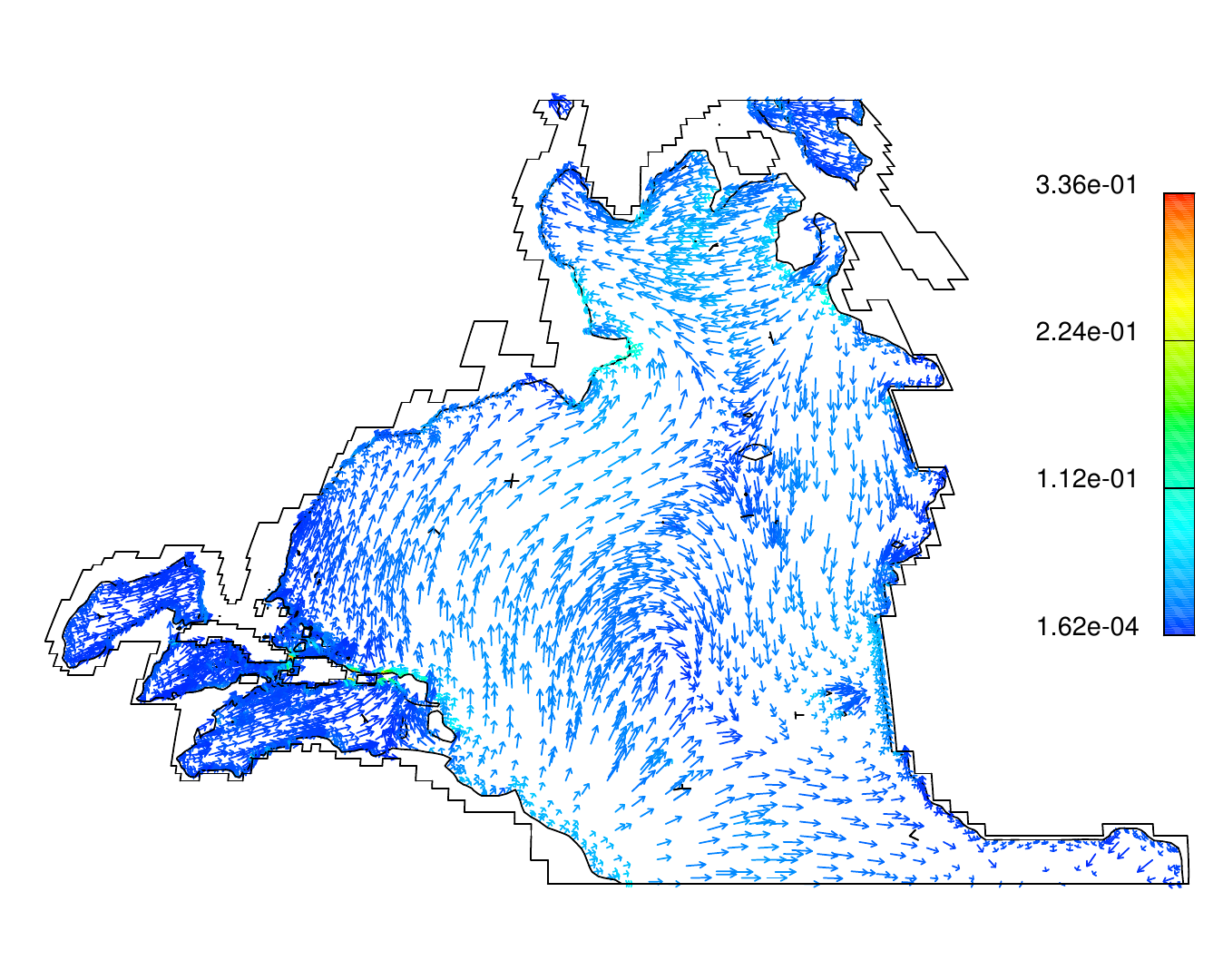}
\end{center}
\caption{Currents at 1500 m depth}
\label{velo1500}
\end{figure}
\begin{figure}[H]
\begin{center}
\includegraphics[width=11cm]{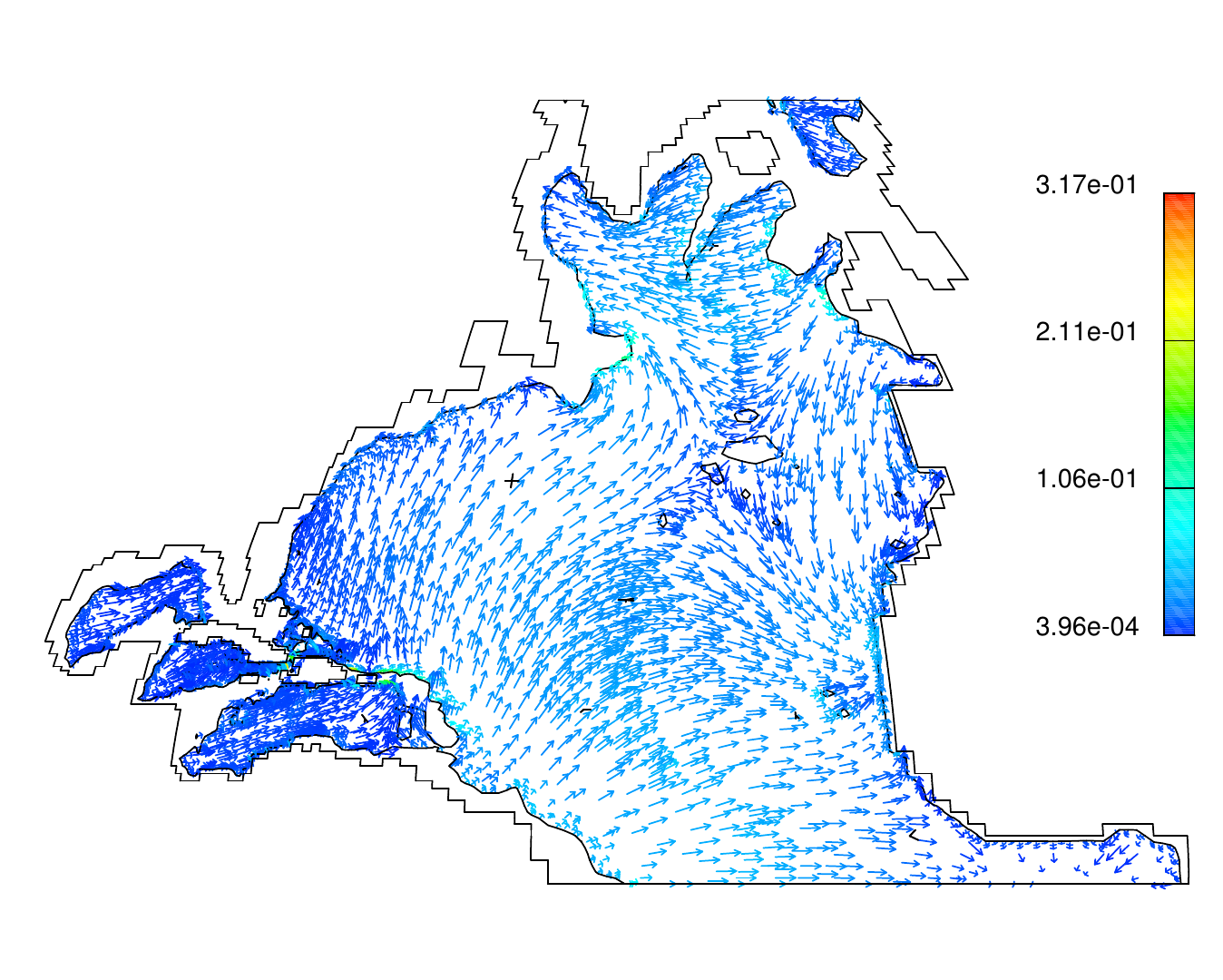}
\end{center}
\caption{Currents at 2000 m depth}
\label{velo2000}
\end{figure}
\begin{figure}[H]
\begin{center}
\includegraphics[width=11cm]{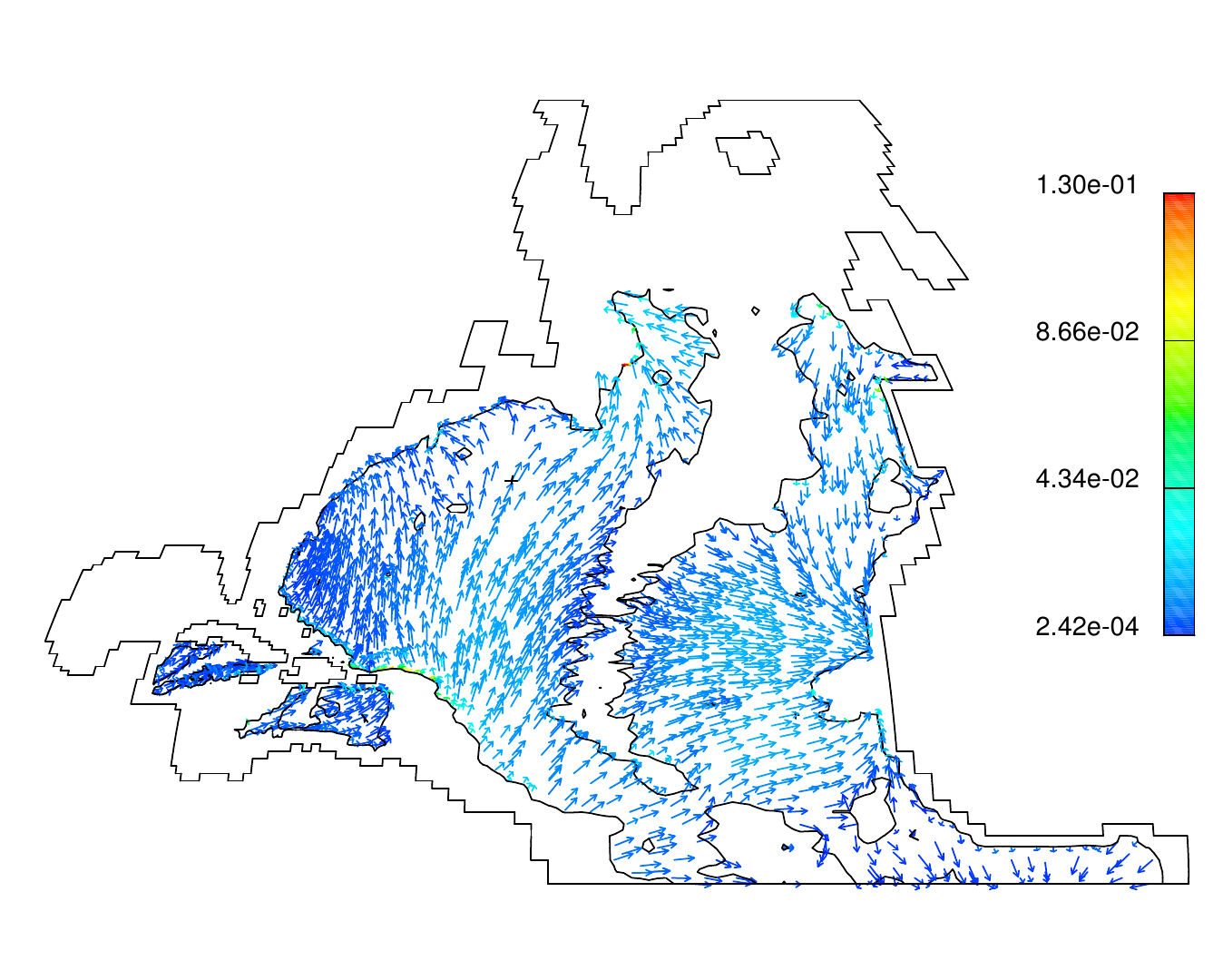}
\end{center}
\caption{Currents at 4000 m depth}
\label{velo4000}
\end{figure}
\begin{figure}[H]
\begin{center}
\includegraphics[width=11cm]{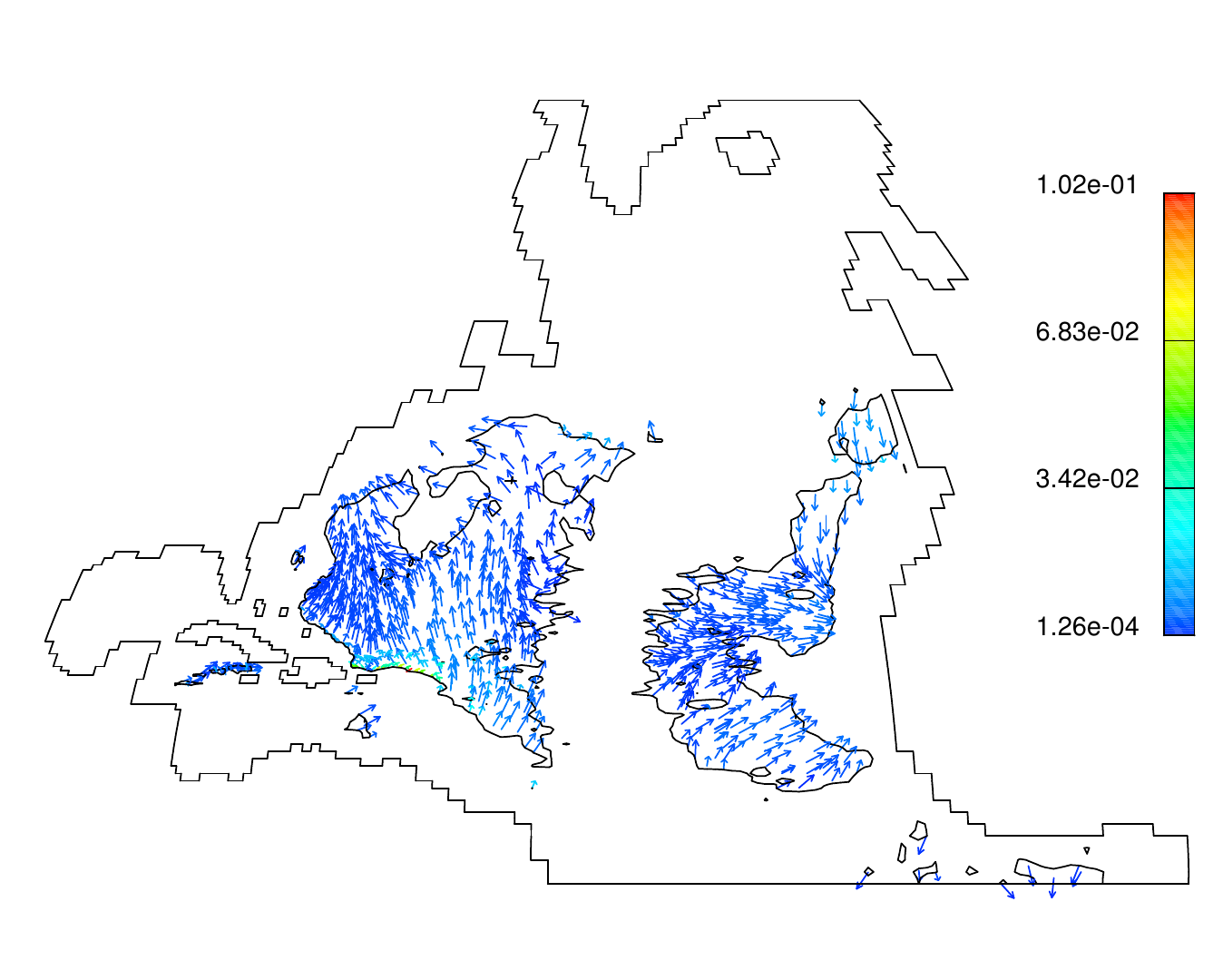}
\end{center}
\caption{Currents at 5000 m depth}
\label{velo5000}
\end{figure}

The next figures show the streamlines associated to the previous velocity field. The upwelling and
downwelling are shown near the American and the African coasts.
\begin{figure}[H]
\begin{center}
\includegraphics[width=15cm]{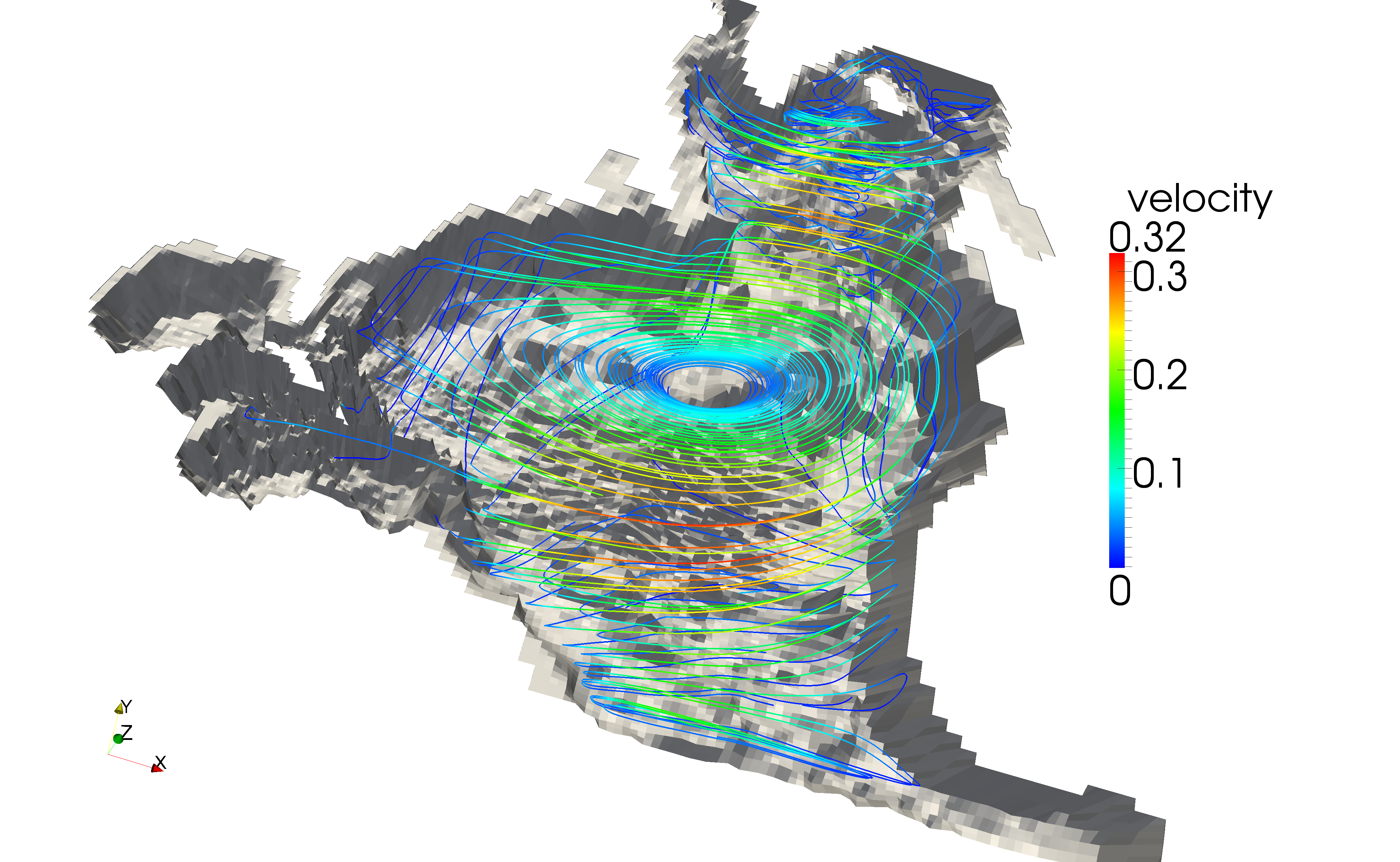}
\end{center}
\caption{Computed streamlines in the North Atlantic}
\label{str0}
\end{figure}
\begin{figure}[H]
\begin{center}
\includegraphics[width=15cm]{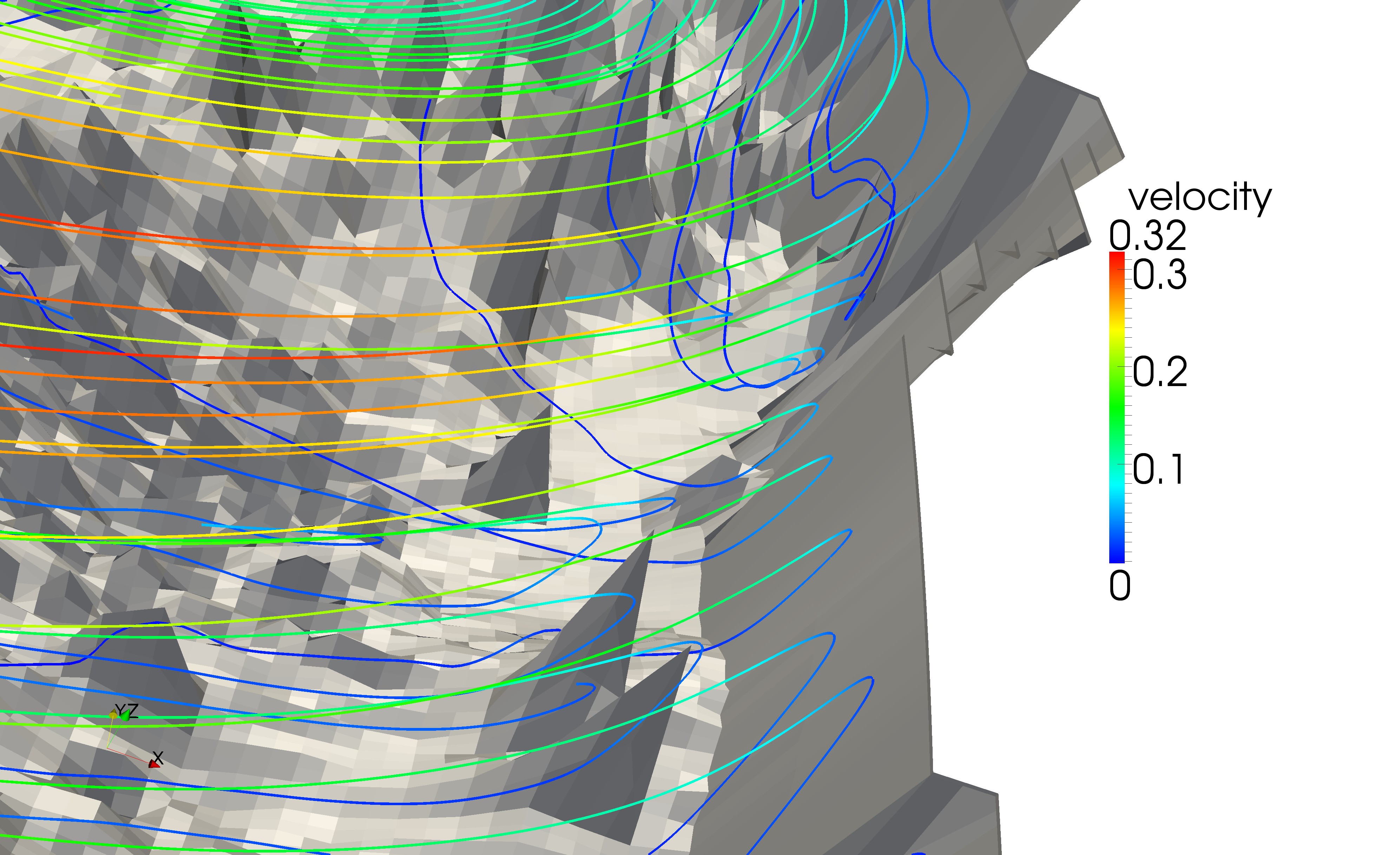}
\end{center}
\caption{Upwelling on the African coast}
\label{str1}
\end{figure}
\begin{figure}[H]
\begin{center}
\includegraphics[width=15cm]{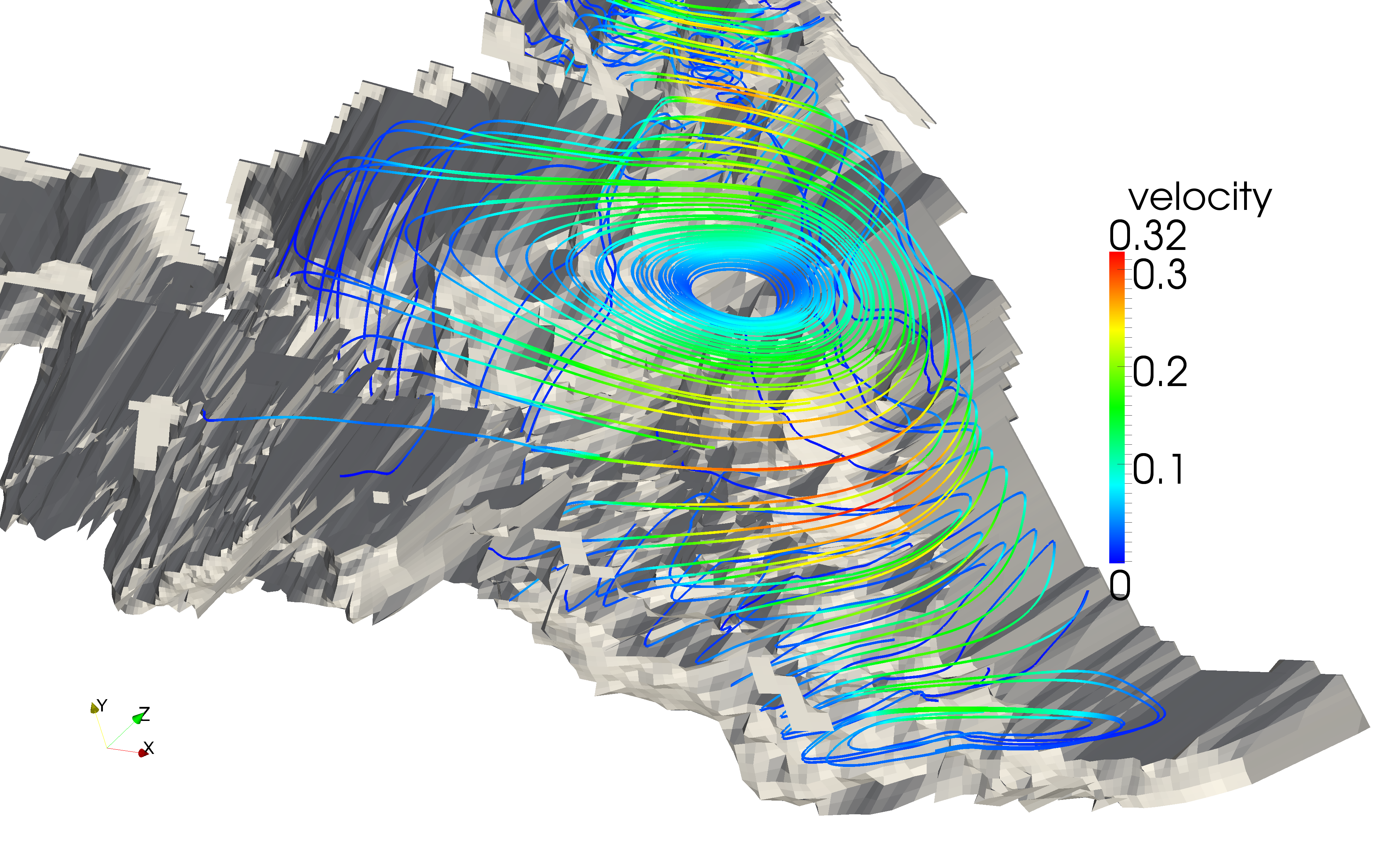}
\end{center}
\caption{Down and upwelling on the American coast}
\label{str2}
\end{figure}
Let us mention that the travel time for a particle of water, at the surface, to cross the Atlantic
from the African coast to the South American coast is about 250 days. The travel time of the same
particle to come back from the South American coast to the African one through the deep water is
about 12 years.

\begin{figure}[H]
\begin{center}
\includegraphics[width=13cm]{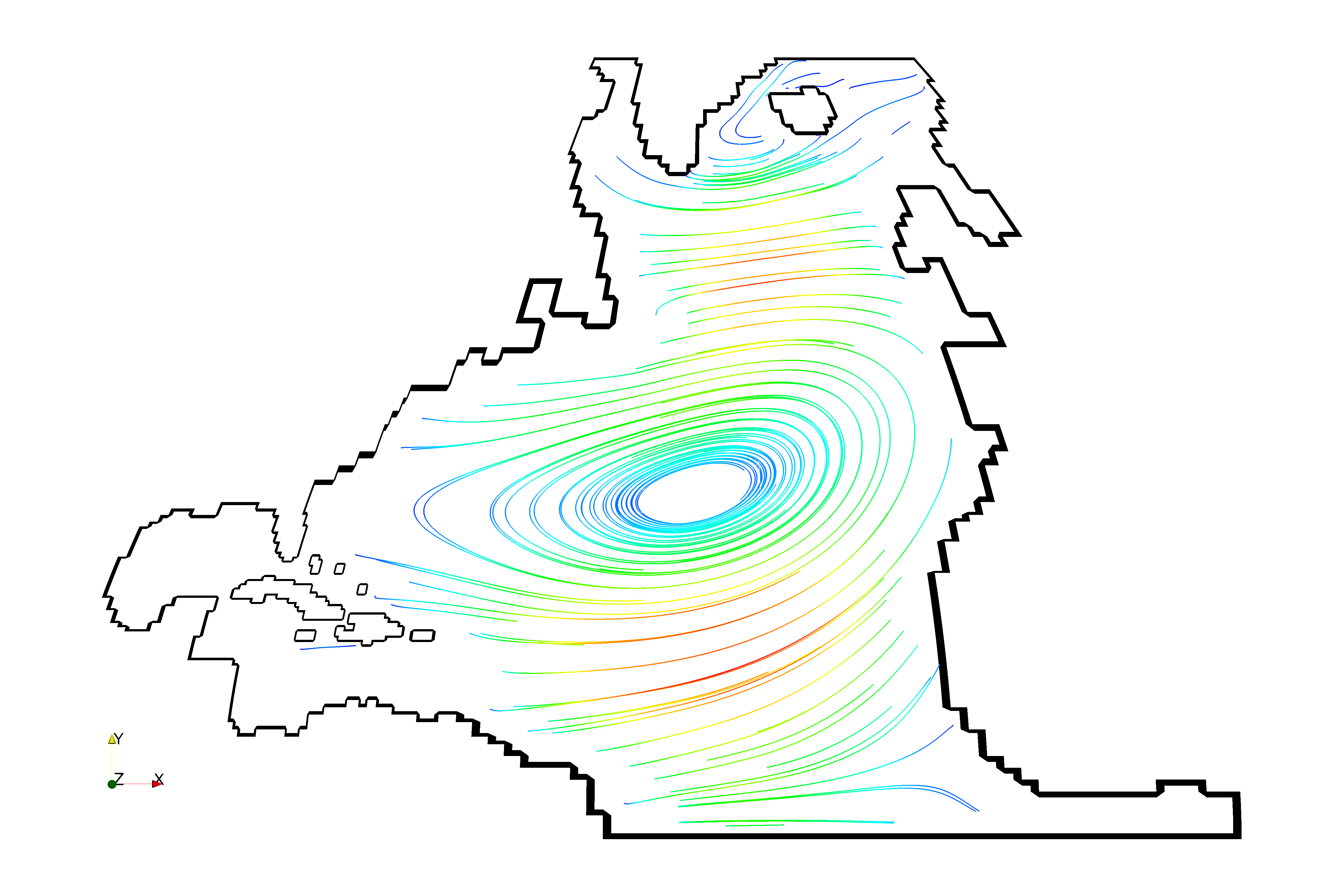}
\end{center}
\caption{Streamlines at the surface}
\label{str_sur}
\end{figure}

\begin{figure}[H]
\begin{center}
\includegraphics[width=12cm]{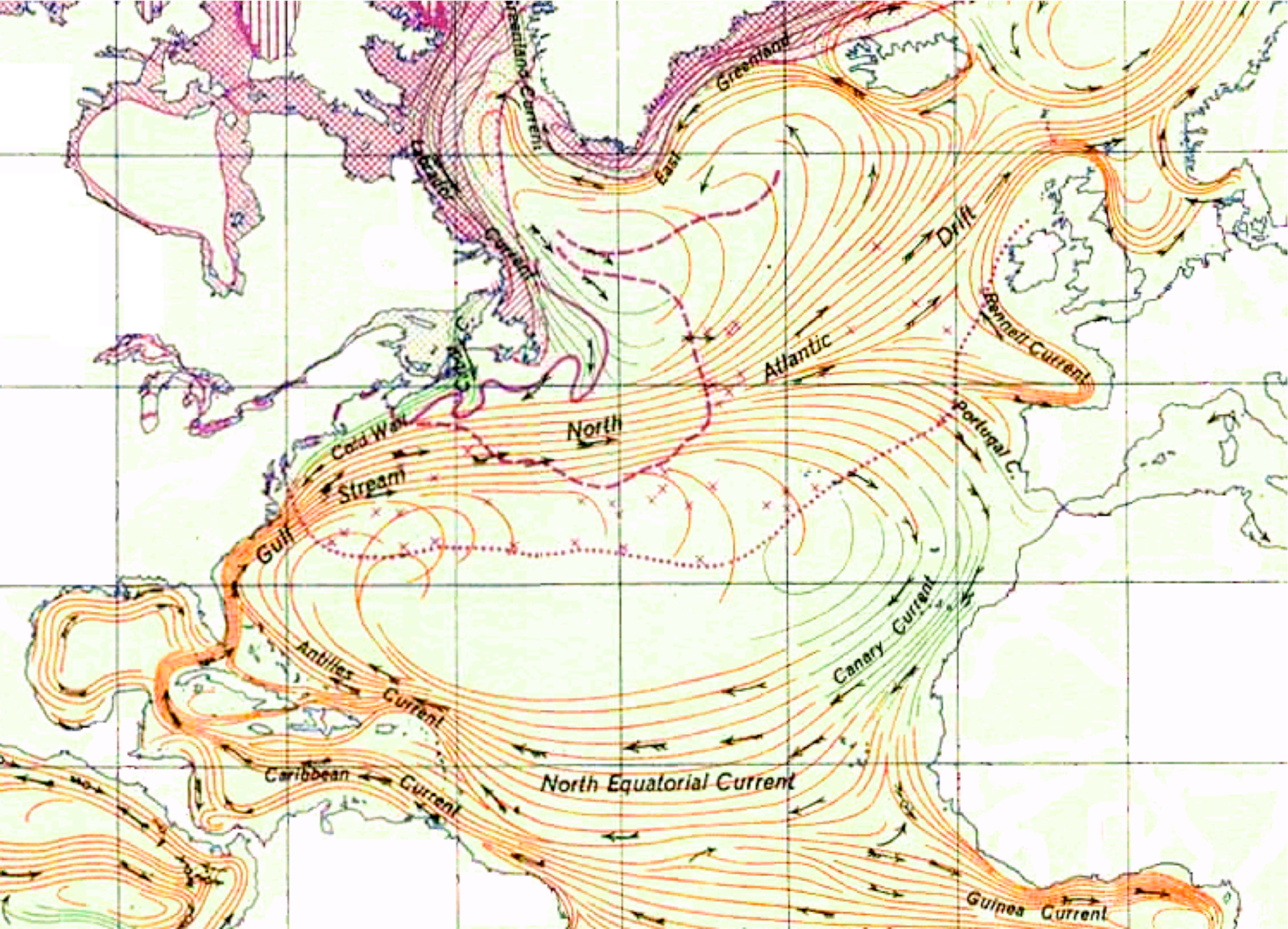}
\end{center}
\caption{North Atlantic Gyre \cite{wiki} }
\label{str_wiki}
\end{figure}

Finally a representation of the surface streamlines obtained with our model is 
represented in figure \ref{str_sur}. This figure can be compared with figure 
\ref{str_wiki} obtained on the web site \cite{wiki}.

The main difference between these two figures are in the region of the Caribbean 
Islands. The reasons for these differences are
\begin{itemize}
 \item Our model do not take into account the temperature and the salinity. 
Therefore we do not have buoyancy effects.
 \item Because of our boundary conditions on the Equator, we could not take into 
account of the South Equatorial current along the East coast of South America. 
For this we could have to compute the currents in the entire Atlantic, 
but we could not have access to the data for this.
\end{itemize}
As already mentioned, our model is able to precisely describe the down, and upwelling near
the coast. Moreover it allows to obtain consistent results with measurements, and demonstrate that
the
winds are the main driving forces for the global dynamic in the oceans. \medskip \\
\textbf{Acknowledgements}. We are grateful to the CSCS  (Swiss National Supercomputing Center,
http://www.cscs.ch) for the use of the computing facilities for a preliminary version of our
software. We are also grateful to Mario Valle at CSCS for his hints in the use of graphics
softwares.
%

%
\end{document}